\documentclass[a4paper,11pt,leqno]{amsart}
\usepackage{amsfonts}
\usepackage{amssymb}
\usepackage{amsmath}
\usepackage[dvipdfmx]{graphicx}
\usepackage{color}
\usepackage{url}
\usepackage{hyperref}

\newtheorem{prop}{Proposition}
\newtheorem{rem}{Remark}

\newcommand\Al{\alpha}
\newcommand\B{\beta}

\newcommand\G{\gamma}

\newcounter{num} 
\newcommand{\Fg}[1][]{\thenum}

\begin{document}



\centerline{}\bigskip

\centerline{}\bigskip

\centerline {\Large{\bf 
Wasan geometry with the division by 0
}}\centerline{}%
\centerline {\Large{\bf }}

\bigskip 

\begin{center}
{\large Hiroshi Okumura}

\centerline{}
\end{center}

\bigskip


\medskip
\textbf{(2010)Mathematics Subject Classification: }01A27, 03C99, 51M04 

\medskip
\textbf{Abstract.} 
Results in Wasan geometry of tangents circles can still be 
considered in a singular case by the division by 0. 
\bigskip

\section{Introduction}

Japanese mathematics in Edo period is called Wasan. 
Wasan geometry considers some relationships which arise when some 
elementary figures such as lines and circles get together. The result 
does not consider the degenerate case explicitly where the circles 
are lines or points. This is a singular case for the parameters 
expressing the circles. In this paper we show that we can still 
consider such a case with the definition of the division by 0. We 
also show that we can consider by manipulating equations of the circles 
with no consideration of limit and the result obtained in the ordinary 
case is still true in the singular case. 

\section{The division by 0}

In this section we give a brief introduction of the division by 0. Let us 
consider the function 
$F: {\rm R}\times {\rm R}\rightarrow {\rm R}$ satisfying 

(i) $F(a,b)=a/b$, if $b\not=0$,  \hskip23mm (ii) $F(a,b)F(c,d)=F(ac,bd)$. 

\noindent
Obviously $F$ is a generalization of the ordinary fraction $a/b$. While 
$F(a,0)=F(a,0)(2/2)=F(2a,0)=F(a,0)(2/1)=2F(a,0)$ gives $F(a,0)=0$. Hence 
we define as follows: 

\medskip
\qquad\qquad (d1) {\it $a/0=0$ for any real number $a$} \cite{kmsy14}.  
\medskip

\noindent
Notice that $F(a,0)=0$ is obtained in any field with characteristic 
different from $2$. 

Any circle or any line has an equation $S(x,y)=a(x^2+y^2)+2gx+2fy+c=0$. 
If $S(x,y)=0$ expresses a circle, its radius is given by 
$$
R=\sqrt{\frac{g^2+f^2-ac}{a^2}}. 
$$
This implies $R=0$, if $a=0$ by (d1). Therefore we define as follows: 

\medskip
\qquad\qquad (d2) {\it If we consider a line as a circle, its radius equals 0 
} \cite{ssdzcm}. 
\medskip

%
%
%

\begin{rem}\label{r1} {\rm 
We get $\tan\pi/2=0$ by (d1), i.e., the slope of a perpendicular equals 
0 if we assume the definition of the division by 0. Therefore we can 
consider that the orthogonality and the tangency are the same 
\cite{MOS17}.}
\end{rem}


\section{The proposition}

Generalizing a problem in Wasan geometry in \cite{aida370}, 
we get the following proposition (see Figure \ref{fprop}). 

\begin{prop}[\cite{OKSJM171}]\label{p2}
Let $\Al$, $\B$, $\G$ be circles of radii $a$, $b$, $c$, respectively.
If $s$ and $t$ are tangents of $\B$ parallel to each other, $\Al$ touches 
$s$ from the same side as $\B$ and $\B$ externally, and $\G$ touches $t$ 
from the same side as $\B$ and $\Al$ and $\B$ externally, then the following 
relation holds:
\begin{equation}\label{eqsth}
c=\frac{b^2}{4a}. 
\end{equation}
\end{prop}

\medskip
\begin{center}
\includegraphics[clip,width=90mm]{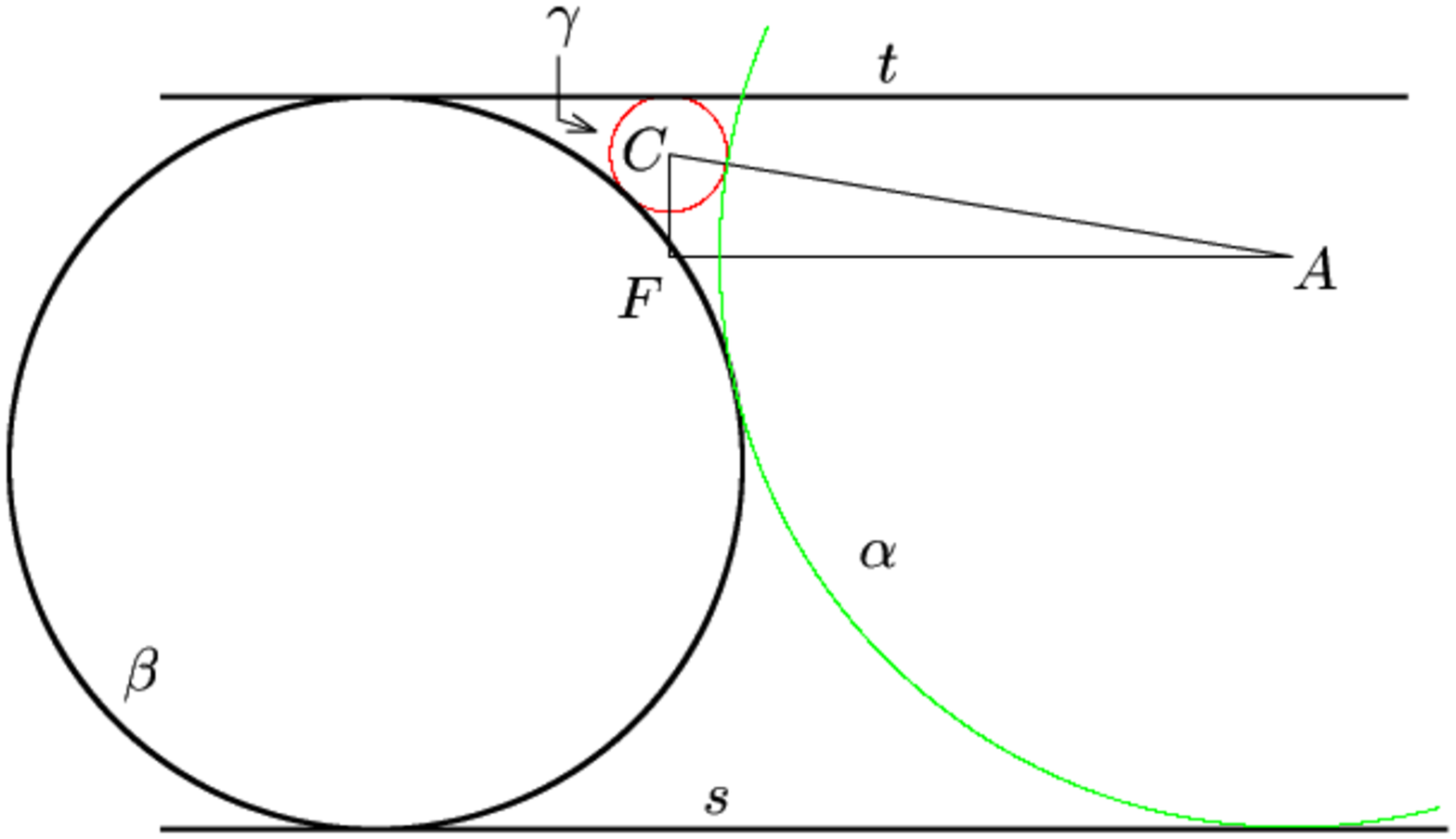}\refstepcounter{num}\label{fprop}\\
\medskip
Figure \Fg .
\end{center}
\medskip

\begin{center}
\includegraphics[clip,width=70mm]{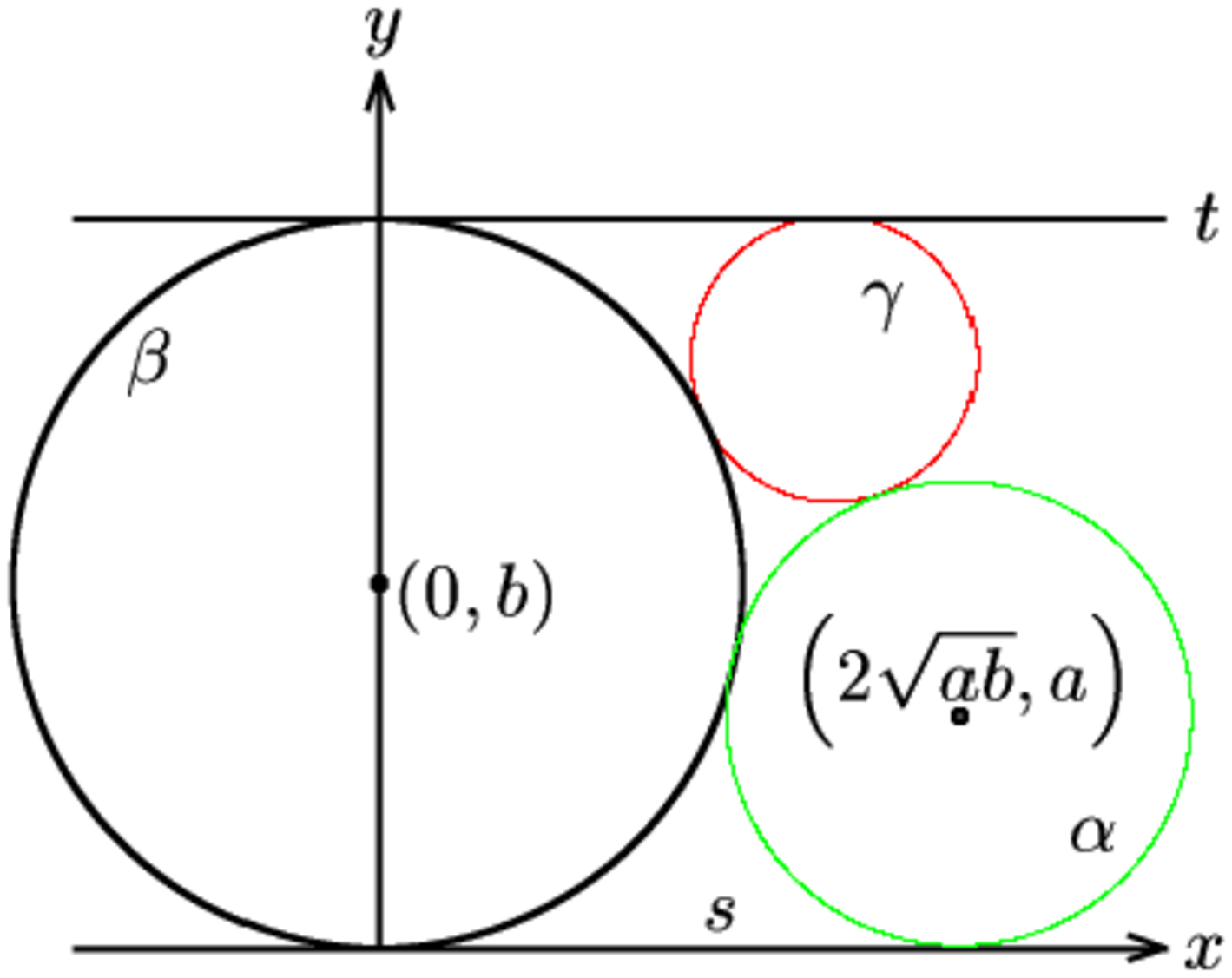}\refstepcounter{num}\label{fax2}\\
\medskip
Figure \Fg .  
\end{center}
\medskip

We now consider the case in which the circle $\Al$ is a point or a line. 
It is equivalent to $a=0$ by (d2). We setup a rectangular coordinate 
system with origin at the point of tangency of the circle $\B$ and the line 
$s$ so that the centers of the circles $\B$ and $\Al$ have coordinates 
$(0,b)$ and $\left(2\sqrt{ab},a\right)$, respectively (see Figure 
\ref{fax2}). Then $\Al$ has an equation 
\begin{equation}\label{eqa1}
\left(x-2\sqrt{ab}\right)^2+(y-a)^2-a^2=0. 
\end{equation}
The equation is arranged as 
\begin{equation}\label{eqa2}
\frac{x^2+y^2}{\sqrt{a}}-4x\sqrt{b}-2\sqrt{a}(y-2b)=0,  
\end{equation}
and 
\begin{equation}\label{eqa3}
\frac{x^2+y^2}{a}-4x\sqrt{\frac{b}{a}}-2(y-2b)=0. 
\end{equation}
If $a=0$, then the equations \eqref{eqa1}, \eqref{eqa2}, \eqref{eqa3} 
imply
\begin{equation}\label{eqax1}\nonumber
x^2+y^2=0, 
\end{equation}
\begin{equation}\label{eqax2}\nonumber
x=0,  
\end{equation}
\begin{equation}\label{eqax3}\nonumber
y=2b, 
\end{equation}
respectively by (d1). 
The last three equations show that $\Al$ is the origin, the $y$-axis, the line 
$t$, respectively. Notice that we can consider that the $y$-axis touches 
the circle $\B$ by Remark \ref{r1}. Therefore the three conclusions are 
reasonable. 

We now consider the circle $\G$ in the same case. 
It has an equation $\left(x-2\sqrt{bc}\right)^2+(y-2b+c)^2=c^2$. 
Since $c=b^2/(4a)$, the equation is arranged as 
\begin{equation}\label{eqg1} \nonumber
a(x^2+(y-2b)^2)-2bx\sqrt{ab}+\frac{b^2y}{2}=0, 
\end{equation} 
\begin{equation}\label{eqg2} \nonumber
\sqrt{a}(x^2+(y-2b)^2)-2bx\sqrt{b}+\frac{b^2y}{2\sqrt{a}}=0, 
\end{equation} 
\begin{equation}\label{eqg3} \nonumber
x^2+(y-2b)^2-2bx\sqrt{\frac{b}{a}}+\frac{b^2y}{2a}=0. 
\end{equation} 
If $a=0$, the three equations give 
\begin{equation}\label{eqgx1} \nonumber
y=0, 
\end{equation} 
\begin{equation}\label{eqgx2} \nonumber
x=0, 
\end{equation} 
\begin{equation}\label{eqgx3} \nonumber
x^2+(y-2b)^2=0, 
\end{equation} 
respectively by (d1). Hence $\G$ is the $x$-axis, the $y$-axis, the point 
$(0,2b)$, respectively.

\medskip
\begin{minipage}{0.31\hsize}
\begin{center}
\vskip2mm
\includegraphics[clip,width=37mm]{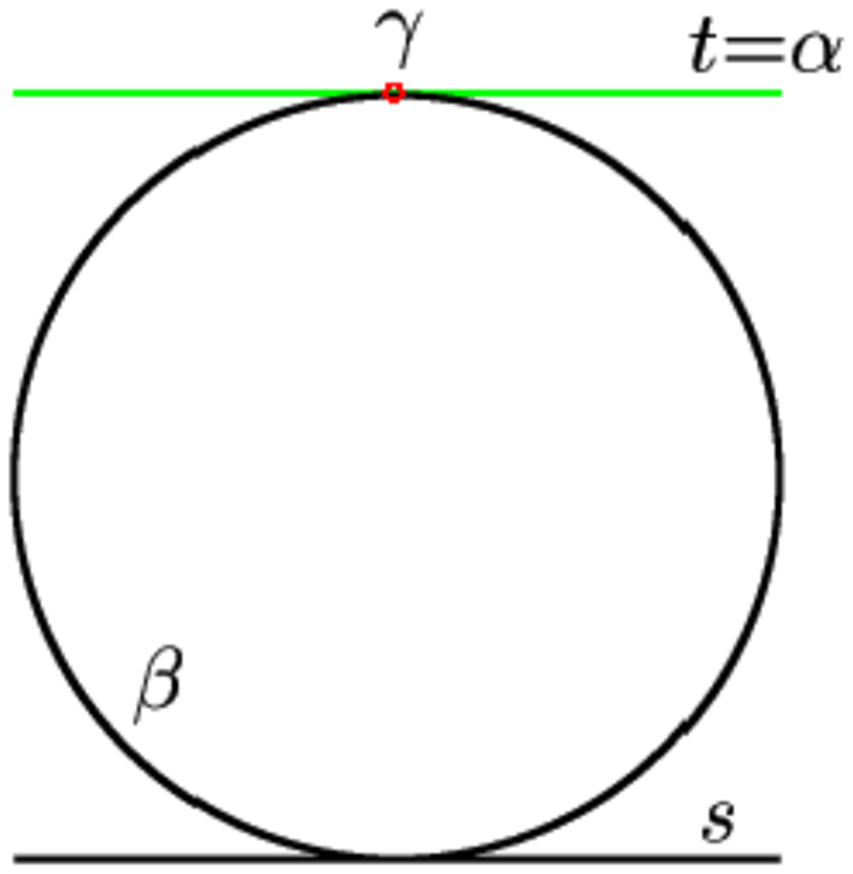}\refstepcounter{num}\label{f3}\\
\vskip2mm
Figure \Fg .
\end{center}
\end{minipage}
\begin{minipage}{0.31\hsize}
\begin{center}
\vskip2mm
\includegraphics[clip,width=41mm]{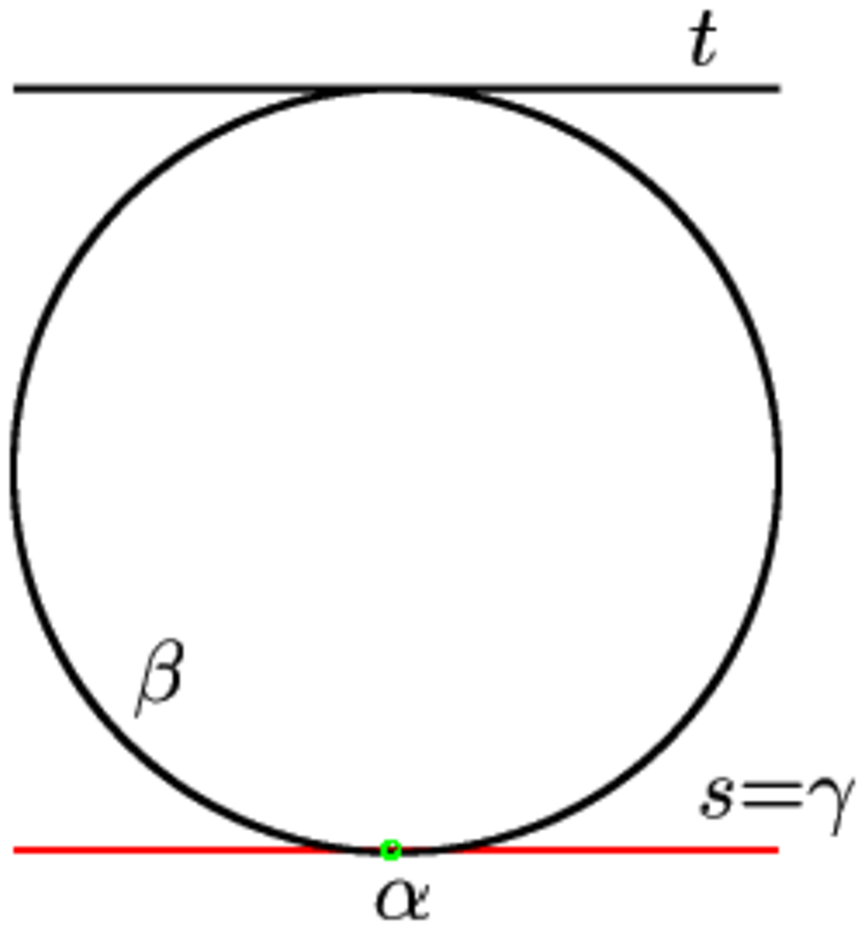}\refstepcounter{num}\label{f4}\\
\vskip2mm
Figure \Fg .
\end{center}
\end{minipage}
\begin{minipage}{0.34\hsize}
\begin{center}
\vskip2mm
\includegraphics[clip,width=35mm]{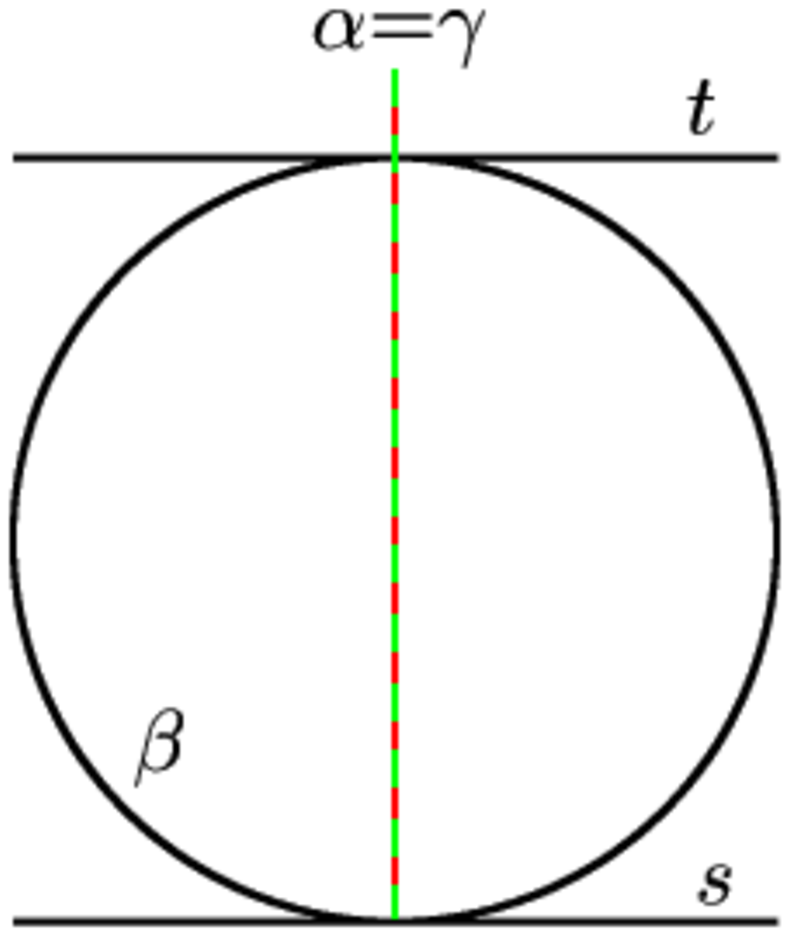}\refstepcounter{num}\label{f5}\\
\vskip2mm
Figure \Fg .
\end{center}
\end{minipage}
\medskip

If $\Al$ approaches to $t$, then $\G$ approaches to the point $(0,2b)$. 
Therefore we can easily consider that $\G$ is $(0,2b)$ if $\Al$ coincides 
with $t$ 
(see Figure \ref{f3}). Symmetrically $\G$ is the line $s$, if $\Al$ 
is the origin (see Figure \ref{f4}). In the rest of the case, both $\Al$ 
and $\G$ coincide with the $y$-axis (see Figure \ref{f5}). 
In all the cases the circle $\G$ is a point or a line, i.e., 
$c=0$ by (d2). Therefore \eqref{eqsth} still holds in all the three cases. 

\section{Conclusion} 

Mathematics is made upon the postulates. Any postulate should be taken into 
consideration, if it gives reasonable conclusions with new insights. 
The three cases, where the circle $\Al$ being a point or a line, can be 
obtained simply and naturally. While one of the three cases in which both 
of the circles $\Al$ and $\G$ coincide with the $y$-axis can not be obtained 
without the definition of the division by 0. Therefore the definition also 
gives us new insights of mathematics.

\bigskip

MAEBASHI GUNMA, 371-0123, JAPAN 

\textit{E-mail address}: \texttt{okmr@protonmail.com}
\end{document}